\documentclass{article}

\usepackage{amsmath,amstext,amsgen,amsbsy,amsopn,amsfonts,graphicx,
theorem,amssymb,shadow,psfrag,subfigure,epsfig,color,cite,comment,url}

\newtheorem{thm}{Theorem}
\newcommand{\floor}[1]{\left\lfloor #1 \right\rfloor}
\newcommand{\R}{\mathbb{R}}

\definecolor{blue}{rgb}{0,0,1}
\definecolor{red}{rgb}{1,0,0}

\begin{document}

\title{The Extraordinary SVD}
\author{Carla D. Martin and Mason A. Porter}  
\date{}
\maketitle

%\footnote{Carla D. Martin is an Associate Professor in the Department of 
%Mathematics and Statistics at James Madison University.}\mbox{ } \footnote{Mason Porter is a University Lecturer in the Mathematical Institute at the University of Oxford.  He is a member of the Oxford Centre for Industrial and Applied Mathematics (OCIAM), a Tutorial Fellow of Somerville College, and a member of Oxford's CABDyN Complexity Centre.

%%%%%%%%%%

%\MAP{is there an arXiv entry for your papers [29,37]?  Update for Ref [4]?  (the readers of the notices don't have access to this)}  
%\CDM{No arXiv entry for those papers.  But I'm guessing they will be published before this
%one gets published, so they can be updated later.}

%\MAP{when I move from bibtex to manual bbl, need to remove "p." from refs [46,54]}

%%%%%%

\begin{abstract}The singular value decomposition (SVD) is a popular matrix factorization that
has been used  widely in applications ever since an efficient algorithm for its computation was developed in the 1970s.  In recent years, the SVD has become even more prominent due to a surge in applications and increased computational memory and speed.  

To illustrate the vitality of the SVD in data analysis, we highlight three of its lesser-known yet fascinating applications: the SVD can be used to characterize political positions of Congressmen, measure the growth rate of crystals in igneous rock, and examine entanglement in quantum computation.  We also discuss higher-dimensional generalizations of the SVD, which have become increasingly crucial with the newfound wealth of multidimensional data and have launched new
research initiatives in both theoretical and applied mathematics.  With its bountiful theory and applications, the SVD is truly extraordinary.
\end{abstract}

%%%%%%

\section{In the Beginning, There is the SVD.}

Let's start with one of our favorite theorems from linear algebra and what is perhaps the most important theorem in this paper.

\begin{thm}Any matrix $A\in\R^{m\times n}$ can be factored into a \emph{singular value decomposition} (SVD),
\begin{equation}
\label{svd}
	A = USV^T\,,
\end{equation}
where $U\in\R^{m \times m}$ and $V\in\R^{n \times n}$ are orthogonal matrices (i.e., $UU^T = VV^T = I$) and $S\in\R^{m\times n}$ is diagonal with $r=\mathrm{rank}(A)$ leading positive diagonal entries.  The 
$p$ diagonal entries of $S$ are usually denoted by $\sigma_i$ for $i=1,\ldots,p$, where $p=\min\{m,n\}$, and $\sigma_i$ are called the \emph{singular values} of $A$.  The singular values are the square roots of the nonzero eigenvalues of both $AA^T$ and $A^TA$, and they satisfy the property
$\sigma_1\geq\sigma_2\geq\dots\geq\sigma_p$.
\end{thm}
See Ref.~\cite{strang} for a proof.

Equation (\ref{svd}) can also be written as a sum of rank-1 matrices,
\begin{equation}
\label{svd2}
	A = \sum_{i=1}^r \sigma_i u_iv_i^T\,,
\end{equation}
where $\sigma_i$ is the $i$th singular value, and $u_i$ and $v_i$ are the $i$th columns of $U$ and $V$.  

Equation \eqref{svd2} is useful when one wants to estimate $A$ using a matrix of lower rank \cite{EY}.  

\begin{thm}\label{thm-EY} (Eckart-Young)
Let the SVD of $A$ be given by \eqref{svd}.  If $k<r=\mathrm{rank}(A)$ and
$\displaystyle A_k = \sum_{i=1}^k \sigma_i u_iv_i^T,$
then 
\begin{align}
	\min_{\mathrm{rank}(B)=k}||A-B||_2=||A-A_k||_2=\sigma_{k+1}\,.
\end{align}
\end{thm}
See Ref.~\cite{vanloan} for a proof.

The SVD was discovered over 100 years ago independently
by Eugenio Beltrami (1835--1899) and Camille Jordan (1838--1921) \cite{Stewart93}.  James Joseph Sylvester (1814--1897), Erhard Schmidt (1876--1959), and Hermann Weyl (1885-1955) also discovered the SVD using different methods \cite{Stewart93}.  The development in the 1960s of practical methods for computing the SVD transformed the field of numerical linear algebra.  One method of particular note is the Golub and Reinsch algorithm from 1970 \cite{GolubReinsch}.  See Ref.~\cite{dhillonsvd} for an overview of properties of the SVD and methods for its computation.  See the documentation for the Linear Algebra Package (LAPACK) \cite{LAPACK} for details on current algorithms to calculate the SVD for dense, structured, or sparse matrices.

Since the 1970s, the SVD has been used in an overwhelming number of applications.  The 
SVD is now a standard topic in many first-year applied mathematics graduate courses and occasionally appears in the undergraduate curriculum.   Theorem \ref{thm-EY} is one of the most important features of the SVD, as it is extremely useful in 
least-squares approximations and principal component analysis (PCA).  During the last decade, the theory, computation, and application of higher-dimensional versions of the SVD (which are based on Theorem \ref{thm-EY}) have also become extremely popular among applications with
multidimensional data.  We include a brief description of a higher-dimensional SVD in this article, and invite you to peruse Ref.~\cite{SIREV} and references therein for additional details. 

We will not attempt in this article to summarize the hundreds of applications that use the SVD, and our discussions and reference list should not be viewed as even remotely comprehensive.  Our goal is to summarize a few examples of recent lesser-known applications of the SVD that we enjoy in order to give a flavor of the diversity and power of the SVD, but there are a myriad of others.  We mention some of these in passing in the next section, and we then focus on examples from Congressional politics, crystallization in igneous rocks, and quantum information theory.  We also discuss generalizations of the SVD before ending with a brief summary.

%Better-known applications are discussed in the references.
%We hope that this article is of interest to mathematics
%students of all levels, mathematicians outside the field of linear algebra, as well as 
%mathematicians whose work directly involves the SVD.

%%%%

\section{It's Raining SVDs (Hallelujah)!}

The SVD constitutes one of science's superheroes in the fight against monstrous data, and it arises in seemingly every scientific discipline.

One finds the SVD in statistics in the guise of ``principal component analysis" (PCA), which entails computing the SVD of a data set after centering the data for each attribute around the mean.  Many other methods of multivariate analysis, such as factor and cluster analysis, have also proven to be invaluable \cite{susanbook}.  The SVD per se has been used in chemical physics to obtain approximate solutions to the coupled-cluster equations, which provide one of the most popular tools used for electronic structure calculations \cite{chemphys}.  Additionally, one applies an SVD when diagonalizing the one-particle reduced density matrix to obtain the natural orbitals (i.e., the singular vectors) and their occupation numbers (i.e., the singular values).  The SVD has also been used in numerous image-processing applications, such as in the calculation of {Eigenfaces} to provide an efficient representation of facial images in face recognition \cite{Turk1,Turk2,eigenface}.  It is also important for theoretical endeavors, such as path-following methods for computing curves of equilibria in dynamical systems \cite{dieci}.  The SVD has also been applied in genomics \cite{genome1,genome2}, textual database searching \cite{text}, robotics \cite{robotics}, financial mathematics \cite{fennpca}, compressed sensing \cite{compress}, and more.  

Computing the SVD is expensive for large matrices, but there are now algorithms that offer significant speed-up (see, for example, Refs.~\cite{Berry92largescale,Larsen98lanczosbidiagonalization}) as well as randomized algorithms to compute the SVD \cite{Rokhlin}.  The SVD is also the basic structure for higher-dimensional factorizations that are SVD-like in nature \cite{SIREV}; this has transformed computational multilinear algebra over the last decade.  

%%%%%%%

\section{Congressmen on a Plane.}

In this section, we use the SVD to discuss voting similarities among politicians.  In this discussion, we summarize work from Refs.~\cite{congshort,conglong}, which utilize the SVD but focus predominantly on other items.

Mark Twain wrote in \emph{Pudd'nhead Wilson's New Calendar} that ``It could probably be shown by facts and figures that there is no distinctly American criminal class except Congress" \cite{twain}. 
There are aspects of this snarky comment that are actually pretty accurate, as much of the detailed work in making United States law is performed by Congressional committees and subcommittees. (This differs markedly from parliamentary democracies such as Great Britain and Canada.)
%, in which the legislative process tends to lie directly in the hands of political parties or is conducted in sessions of the entire parliament.  

There are many ways to characterize the political positions of Congressmen.  %For example, one can tabulate individual voting records on selected key issues using, say, interest group ratings.  However, this is subjective by nature and is also not very mathematically appealing.  
An objective approach is to apply data-mining techniques such as the SVD (or other ``multidimensional scaling" methods) on matrices determined by the Congressional roll call.   Such ideas have been used successfully for decades by political scientists such as Keith Poole of UC San Diego and Howard Rosenthal of Princeton University
%, who have made pioneering contributions in this area 
\cite{pr97,voteview}.  One question to ask, though, is what observations can be made using just the SVD.

In Refs.~\cite{congshort,conglong}, the SVD was employed to investigate the ideologies of Members of Congress.  
%The SVD calculations in those papers were motivated by a then-recent paper in which Larry Sirovich performed a similar investigation using United States Supreme Court cases \cite{sirovich}.  Sirovich noted the excellent accuracy of SVDs in reconstructing votes and conducted a fascinating comparison of cases from the ``second Rehnquist Court" (which began August 8, 1994 and lasted over eight years) with those from the 1959--1961 and 1967--1969 Warren Courts.  He also used Shannon information\footnote{If the set $\{p_k\}$ gives the probabilities of each of the possible outcomes, then the entropy (i.e., the mean information) conveyed by a Supreme Court decision is $I = -\sum_k p_k \log_2 p_k$ \cite{sirovich}.} to show that the second Rehnquist Court behaved as if it were composed of 4.68 ideal justices rather than the 9 nonideal ones that actually comprised it.  One can get very far just by considering the usual SVD (rather than more intricate data-mining techniques, some of which take explicit advantage of political information \cite{voteview}), as it provides a simple method of grouping like-minded politicians based on voting data.  
Consider each two-year Congress as a separate data set and also treat the Senate and House of Representatives separately.  Define an $m\times n$ voting matrix $A$ with one row for each of the $m$ legislators and one column for each of the $n$ bills on which legislators voted.  The element $A_{ij}$ has the value $+1$ if legislator~$i$ voted ``yea'' on bill~$j$ and $-1$ if he or she voted ``nay.''  The sign of a matrix element has no bearing {\it a priori} on conservativism versus liberalism, as the vote in question depends on the specific bill under consideration.  %(However, one would expect the votes on the majority of bills to lean in the direction of the majority party.)  
If a legislator did not vote because of absence or abstention, the corresponding element is~$0$.  %Abstentions and absences are treated as the same.  
Additionally, a small number of false zero entries result from resignations and midterm replacements.

Taking the SVD of $A$ allows one to identify Congressmen who voted the same way on many bills.  Suppose the SVD of $A$ is given by 
\eqref{svd2}.  The grouping that has the largest mean-square overlap with the actual groups voting for or against each bill is given by the first left singular vector $u_1$ of the matrix,
the next largest by the second left singular vector $u_{2}$, and so on.  
Truncating $A$ by keeping only the first $k \leq r$ nonzero singular values gives the approximate voting matrix
\begin{equation}
	A_k = \sum_{i=1}^k \sigma_i u_{i} v_{i}^T \approx A\,. \label{full}
\end{equation}
This is a ``$k$-mode truncation" (or ``$k$-mode projection") of the matrix $A$.  By Theorem \ref{thm-EY}, \eqref{full} is 
a good approximation as long as the singular values decay sufficiently rapidly with increasing $i$. 

%The sum of the squares of the errors in the elements in (\ref{full}) is equal to $\displaystyle\sum_{i=k+1}^r \sigma_i^2$, which vanishes in the limit $k \to r$.  This approximation is a good one as long as the singular values decay sufficiently rapidly with increasing $i$.  The $\ell$th term in the singular value decomposition~(\ref{svd2}) accounts for a fraction $\sigma_\ell^2/\sum_{i=1}^r \sigma_i^2$ of the sum of the squares of the elements in the voting matrix.

\begin{figure}
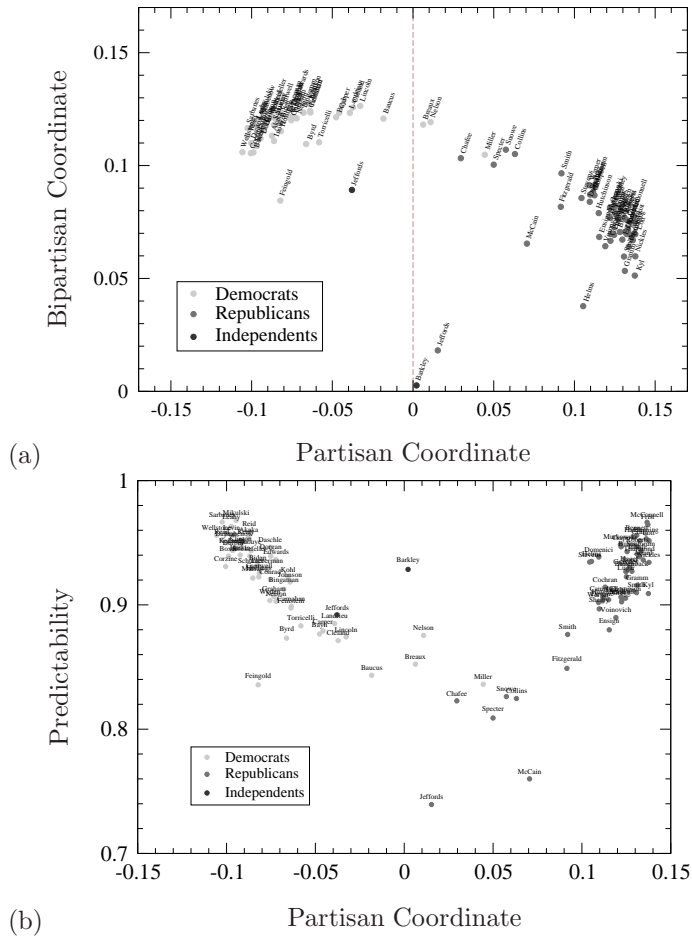

  {\psfragscanon
  \psfrag{107th Senate}{\hspace*{1.2in}107th Senate}
  \psfrag{partisan}[t][t]{Partisan Coordinate}
  \psfrag{partisan coordinate}[t][t]{Partisan Coordinate}
  \psfrag{bipartisan}[b][b]{Bipartisan Coordinate}
  \psfrag{predictability}[b][b]{Predictability}
  %\centerline{
  (a)
\includegraphics[width=0.7\textwidth]{senate107labeledrecolor.eps}\\
  (b) \includegraphics[width=0.7\textwidth]{predict107recolor.eps}}
  %}
\caption{Singular value decomposition (SVD) of the Senate voting
  record from the 107th U.S.~Congress (2001--2002).  (a) Two-mode truncation $A_2$ of
  the voting matrix $A$.  Each point represents a projection of a single
  Representative's votes onto the leading two eigenvectors (labeled
  ``partisan'' and ``bipartisan,'' as explained in the text).  Democrats
  (light dots) appear on the left and Republicans (medium dots) are on the right.
  The two Independents are shown using dark dots .  (b) ``Predictability" of votes cast by Senators in the 107th Congress based on a two-mode truncation of the SVD.  Individual Senators range from 74\% predictable to 97\% predictable.  These figures are modified versions of figures that appeared in Ref.~\cite{conglong}.
}\label{senate}
\end{figure}

A Congressman's voting record can be characterized by just two coordinates \cite{conglong,congshort}, so the two-mode truncation $A_2$ is an excellent approximation to $A$.  One of the two directions (the ``partisan'' coordinate) correlates well with party affiliation for members of the two major parties.  The other direction (the ``bipartisan'' coordinate) correlates well with how often a Congressman votes with the majority.% (i.e., how well he/she plays with others).
\footnote{For most Congresses, it suffices to use a two-mode truncation.  For a few, it is desirable to keep a third singular vector, which can be used to try to encapsulate a North-South divide \cite{conglong,voteview}.}  We show the coordinates along these first two singular vectors for the 107th Senate (2001--2002) in Fig.~\ref{senate}a.  As expected, Democrats (on the left) are grouped together and are almost completely separated from Republicans (on the right).\footnote{Strictly speaking, the partisanship singular vector is determined up to a sign, which is then chosen to yield the usual Left/Right convention.}  The few instances of party misidentification are unsurprising; Conservative Democrats such as Zell Miller [D-GA] appear farther to the right than some moderate Republicans~\cite{org}.  Senator James Jeffords [I-VT], who left the Republican party to become an Independent early in the 107th Congress, appears closer to the Democratic group than the Republican one and to the left of several of the more conservative Democrats.\footnote{Jeffords appears twice in Fig.~\ref{senate}a---once each for votes cast under his two different affiliations.}  
%The behavior of legislators in the House of Representatives can also be described extremely well using two-mode truncations of the SVDs of the voting matrices \cite{conglong,congshort}. 

Equation \eqref{full} can also be used to construct an approximation to the votes %(as opposed to the voters) 
in the full roll call.  Again using $A_2$, one assigns ``yea'' or ``nay'' votes to Congressmen based on the signs of the matrix elements.  Figure \ref{senate}b shows the fraction of actual votes correctly reconstructed using this approximation.  Looking at whose votes are easier to reconstruct gives a measure of the ``predictability'' of the Senators in the 107th Congress.  Unsurprisingly, moderate Senators are less predictable than hard-liners for both parties.  Indeed, the two-mode truncation correctly reconstructs the votes of some hard-line Senators for as many as 97\% of the votes that they cast.

\begin{figure}
{\psfragscanon
  \psfrag{partisan coordinate}[t][t]{Partisan Coordinate}
  \psfrag{bipartisan coordinate}[b][b]{Bipartisan Coordinate}
\centerline{
  \includegraphics[width=.9\textwidth]{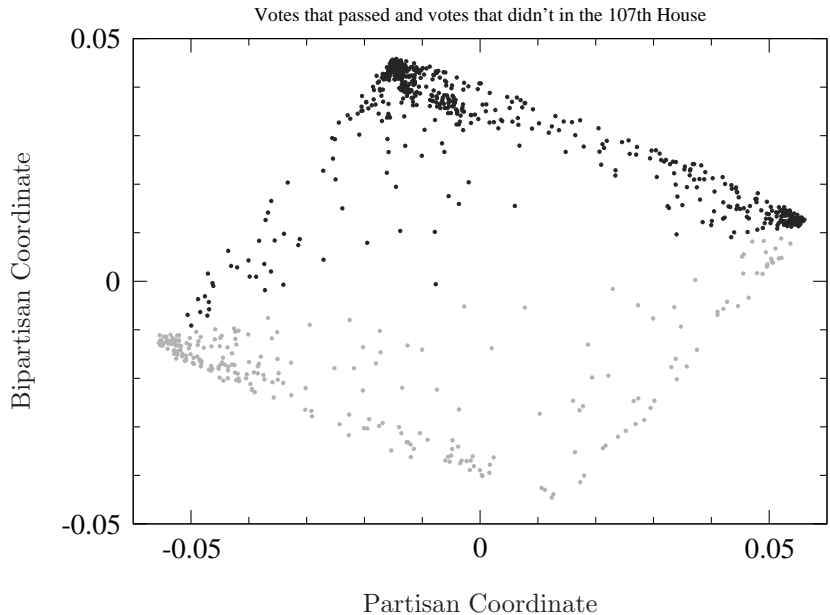}}}
\caption{SVD of the roll call of the 107th House of Representatives projected onto the voting coordinates.  
  %Points represent projections of the votes cast on a
  %measure onto eigenvectors associated with the leading two singular
  %values. 
  There is a clear separation between bills that passed (dark dots) and those that did not (light dots).  The four corners of the plot are interpreted as follows: bills with broad bipartisan support (north) all passed; those supported mostly by the Right (east) passed because the Republicans were the majority party; bills supported by the Left (west) failed because of the Democratic minority; and the (obviously) very few bills supported by almost nobody (south) also failed.  This figure is a modified version of a figure that appeared in Ref.~\cite{conglong}.}
\label{votesvd}
\end{figure}

To measure the reproducibility of individual votes and outcomes, the SVD can be used to calculate the positions of the votes along the partisanship and bipartisanship coordinates (see Fig.~\ref{votesvd}).  One obtains a score for each vote by reconstituting the voting matrix as before using the two-mode truncation $A_2$ and summing the elements of the approximate voting matrix over all legislators.  Making a simple assignment of ``pass'' to those votes that have a positive score and ``fail'' to all others successfully reconstructs the outcome of $984$ of the $990$ total votes (about 99.4\%) in the 107th House of Representatives.  A total of $735$ bills passed, so simply guessing that every vote passed would be considerably less effective.  This way of counting the success in reconstructing the outcomes of votes is the most optimistic one.  Ignoring the values from known absences and abstentions, $975$ of the $990$ outcomes are still identified correctly.  Even the most conservative measure of the reconstruction success rate---in which one ignores values associated with abstentions and absences, assigns individual yeas or nays according to the signs of the elements of $A_2$, and then observes which outcome has a majority in the resulting
roll call---identifies $939$ (about 94.8\%) of the outcomes correctly.  The success rates for other recent Houses are similar \cite{conglong}.

To conclude this section, we remark that it seems to be underappreciated that many political scientists are extremely sophisticated in their use of mathematical and statistical tools.  Although the calculations that we discussed above are heuristic ones, several mathematicians and statisticians have put a lot of effort into using mathematically rigorous methods to study problems in political science.  For example, Donald Saari has done a tremendous amount of work on voting methods \cite{saari}, and (closer to the theme of this article) rigorous arguments from multidimensional scaling have recently been used to study roll-call voting in the House of Representatives \cite{sharad}.

%%%%%%%%%%%%%%%%%%%%%%%%
\section{The SVD is Marvelous for Crystals.}

Igneous rock is formed by the cooling and crystallization of magma.  One interesting aspect of the formation of igneous rock is that the microstructure of the rock is composed of interlocking crystals of irregular shapes.  The microstructure contains a plethora of quantitative information about the crystallization of deep crust---including the nucleation and growth rate of crystals.  In particular, the three-dimensional (3D) \emph{crystal size distribution} (CSD) provides a key piece of information in the study of crystallization rates.  CSD can be used, for example, to determine the ratio of nucleation rate to growth rate.  Both rates are slow in the deep crust, but the growth rate dominates the nucleation rate.  This results in a microstructure composed of large crystals.  See Ref.~\cite{Amenta4} for more detail on measuring growth rates of crystals and Refs.~\cite{marsh98,higgins06} for more detail on this application of the SVD.
 
\begin{figure}
\centering
\mbox{\subfigure[Tetragonal Prism (1:1:5)]{\phantom{XXXXX}\epsfig{file=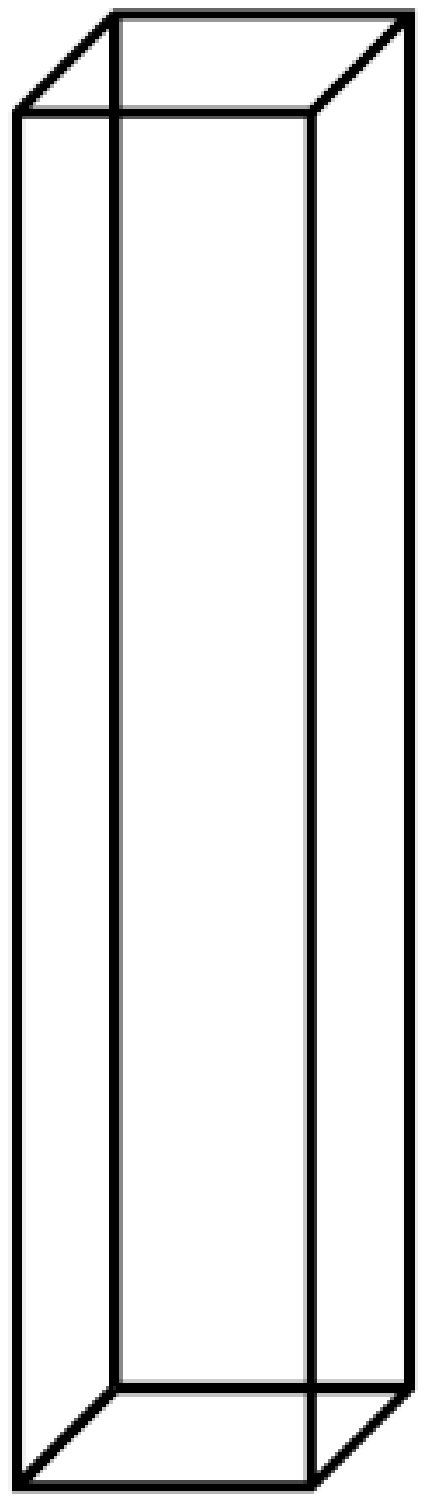,height=2in}\phantom{XXXX}}\quad
\subfigure[Tetragonal Plate (1:5:5)]{\phantom{XXXXX}\epsfig{file=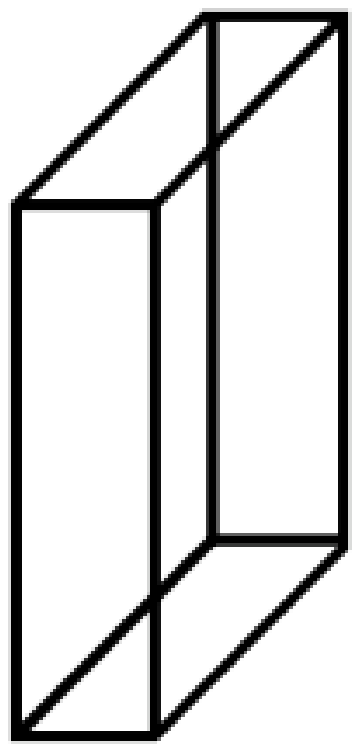,height=2in}\phantom{XXXX}}\subfigure[Orthorhombic cuboid (1:3:5)]{\epsfig{file=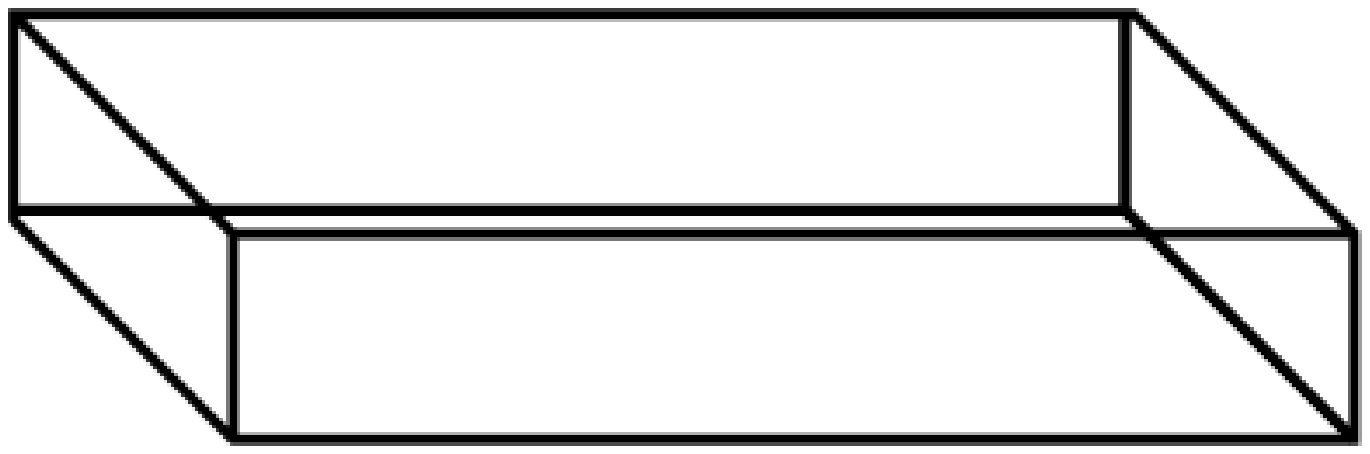,width=2in}}}
\caption{Crystalline structures used to measure grain sizes.  We give the relative sizes of their dimensions in parentheses.}
\label{fig:crystalshapes}
\end{figure}

As the crystals in a microstructure become larger, they compete for growth space and their grain shapes become irregular.  This makes it difficult to measure grain sizes accurately.  CSD analysis of rocks is currently done in two stages.  First,  one takes hand measurements of grain sizes in 2D slices and then computes statistical and stereological corrections to the measurements in order to estimate the actual 3D CSD.  However, a novel recent approach allows one to use the SVD to automatically and directly measure 3D grain sizes that are derived from three specific crystal shapes (prism, plate, and cuboid; see Fig.~\ref{fig:crystalshapes}) \cite{Amenta1}.  Ongoing research involves extending such analysis to more complex and irregular shapes.  Application to real rock microstructures awaits progress in high energy X-ray tomography, as this will allow improved resolution of grain shapes.

The grain sizes are determined by generating databases of 
microstructures with irregular grain shapes in order to compare the estimated CSD of the actual grains to the computed or ideal CSD predicted by the governing equations.  Because the CSDs in many igneous rocks are close to linear \cite{Amenta2,Amenta1}, the problem can be simplified by using governing equations that generate linear CSDs with the following two rate laws.
\begin{enumerate}
\item{\emph{Nucleation Rate Law}: $N(t)=e^{\alpha t}$, where $N$ is the number of new nuclei formed at each time step $t$ and $\alpha$ is the nucleation constant.
}
\item{\emph{Crystal Growth Rate Law}: $G=\Delta L/\Delta t$, where $\Delta L/\Delta t$ is the 
rate of change of a grain diameter per time step.  Grain sizes can be represented by short,
intermediate, or long diameters.  Such diameter classification depends on the relationship
between the rate of grain nucleation and the rate of grain growth.
}
\end{enumerate}

\begin{figure}
\centering
\mbox{\subfigure[Enclosing Ellipsoid]{\epsfig{file=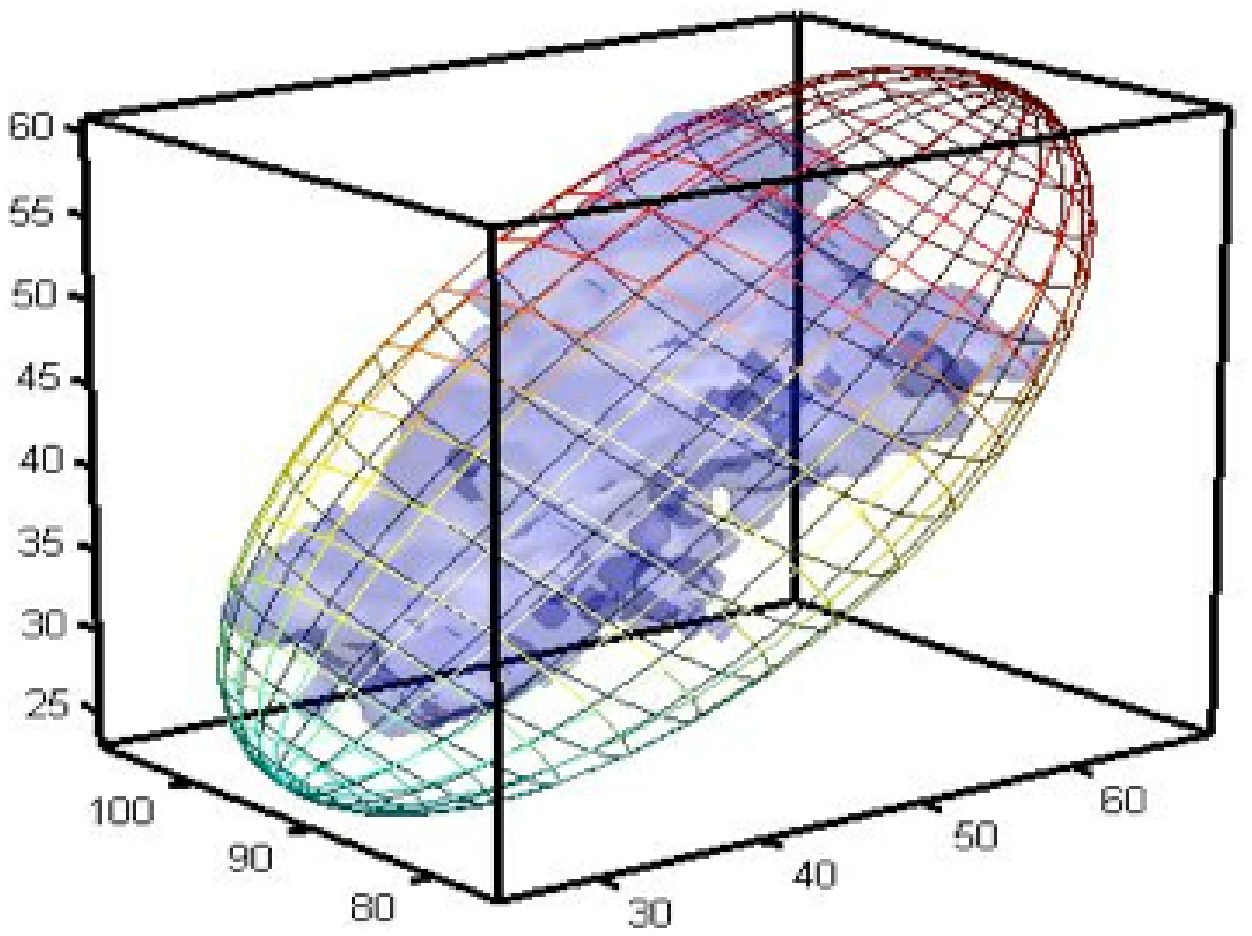,width=2.5in}}\quad
\subfigure[Inscribed Ellipsoid]{\epsfig{file=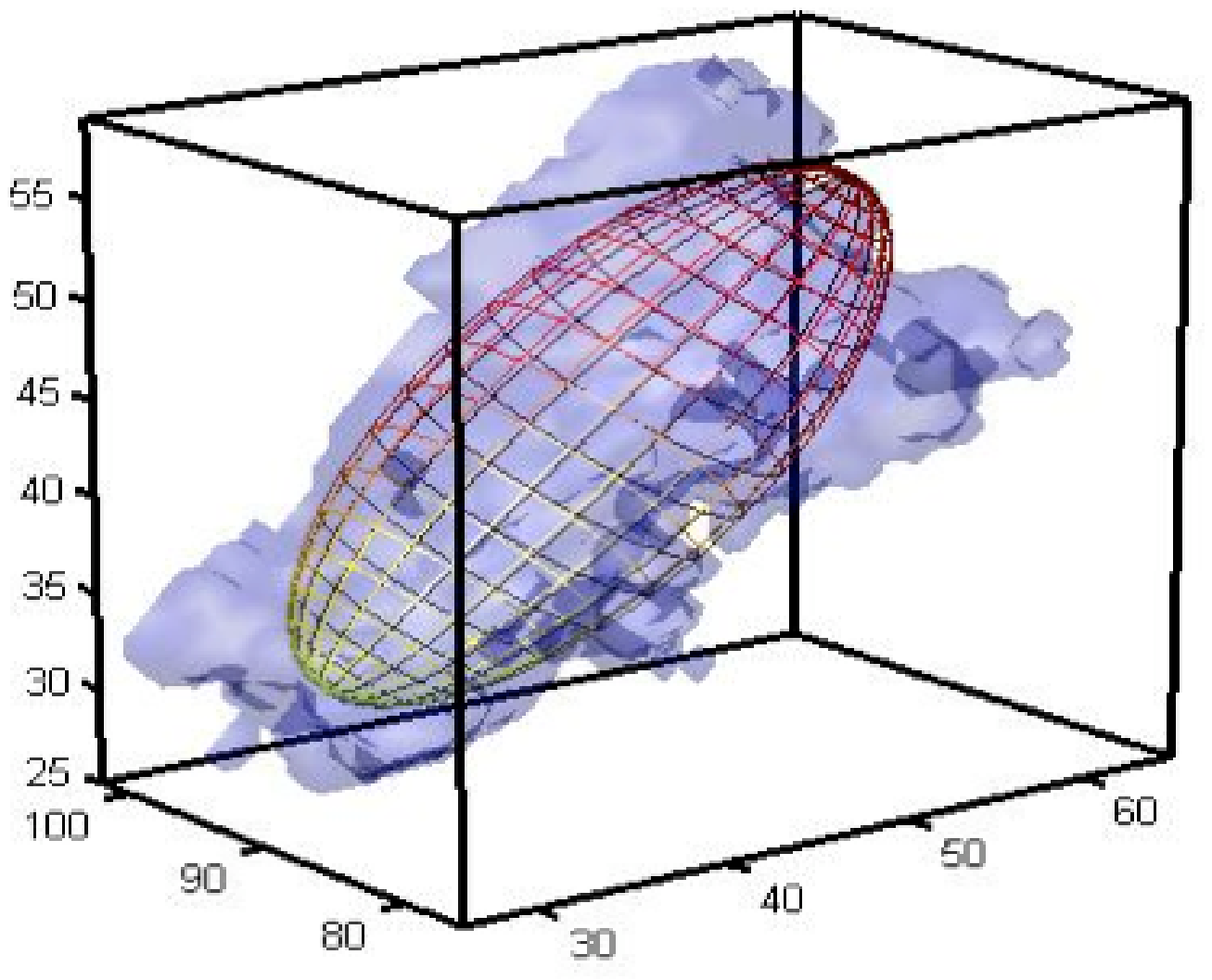,width=2.5in}}}
\caption{Two possible ellipsoids used to approximate grain sizes. Because grain
shapes are irregular, all ellipsoids are triaxial with three unequal diameters. 
}
\label{fig:ellipsoids}
\end{figure}

One uses an ellipsoid to approximate the size and shape of each grain.  There are multiple subjective choices for such ellipsoids that depend on the amount (i.e., the number of
points) of the grain to be enclosed by the ellipsoid.  To circumvent this subjectivity, it is desirable to compare the results of three types of
ellipsoids: the ellipsoid that encloses the entire grain, the ellipsoid that is inscribed within the grain, and the mean of the enclosed and inscribed ellipsoids. See Fig.~\ref{fig:ellipsoids} for an illustration of an enclosing and an inscribed ellipsoid.

The SVD is used in the determination of each of the three types of ellipsoids.
Comparing the CSDs obtained using each of the three types of
ellipsoids with those predicted by the governing equations reveals that the inscribed ellipsoids give the best results.  In particular, one can use an algorithm developed by Nima Moshtagh \cite{moshtagh} that employs the Khachiyan Algorithm \cite{khachiyan} along with the SVD to obtain an ellipsoid that encloses an arbitrary number of points (which is defined by the user).Ê Leonid Khachiyan introduced the ellipsoid method in 1979, and this was the first algorithm for linear programming that runs in polynomial time in the worst case.  Given a matrix of data points $P$ containing a discretized set of 3D points representing the crystal,
one solves
\begin{align}
	\min_{A,c}\log\{\det(A)\} \quad \mbox{subject to} \quad (P_i-c)^TA(P_i-c)\leq 1\,,
\end{align}
where $P_i$ is the $i$th column of $P$, the matrix $A$ contains information about the
shape of the ellipsoid, and $c$ is the center of the ellipsoid.

Note that $P$ in this case is dense, it has size $n\times 3$, and $n\approx 5000$.
Once $A$ and $c$ have been determined, one calculates the $i$th radius of the $D$-dimensional ellipse from the SVD of $A$ using
\begin{align}
	r_i&=1/\sqrt{\sigma_i}\,,
\end{align}
where $\sigma_i$ ($i=1,\ldots, D$) is the $i$th singular value of $A$.  If the SVD of $A$ is given by equation \eqref{svd}, then the orientation of the ellipsoid is given by the rotation matrix $V$.

The major difficulty in such studies of igneous rock is that grain shapes and sizes are irregular due to competition for growth space among crystals.  In particular, they are not of the ideal sizes and shapes that are assumed by crystallization theory.  For example, crystals might start to grow with definite diameter ratios (yielding, for example, the prism, plate, or cuboid in Fig.~\ref{fig:crystalshapes}) but eventually develop irregular outlines.  Current studies \cite{Amenta1} suggest that one of the diameters or radii of the inscribed ellipsoid (as determined from the SVD) can be used as a measure of grain size for the investigation of crystal size distributions, but the problem remains open.  

%%%%%%%%%%%%%
%%%%%%

\section{Quantum Information Society.}

From a physical perspective, information is encoded in the state of a physical system, and a computation is carried out on a physically realizable device \cite{preskill}.  \emph{Quantum information} refers to information that is held in the state of a quantum system.  Research in quantum computation and quantum information theory has helped lead to a revival of interest in linear algebra by physicists.  In these studies, the SVD (especially in the form of the Schmidt decomposition) have been crucial for gaining a better understanding of fundamental quantum-mechanical notions such as entanglement and measurement.

{\it Entanglement} is a quantum form of correlation that is much stronger than classical correlation, and quantum information scientists use entanglement as a basic resource in the design of quantum algorithms \cite{preskill}.  The potential power of quantum computation relies predominantly on the inseparability of multipartite quantum states, and the extent of such interlocking can be measured using entanglement.  

We include only a brief discussion in the present article, but one can go much farther \cite{preskill,pasyou,sch}.  Whenever there are two distinguishable particles, one can fully characterize inseparable quantum correlations using what is known as a ``single-particle reduced density matrix" (see the definition below), and the SVD is crucial for demonstrating that this is the case.  See Refs.~\cite{preskill,pasyou,sch} for lots of details and all of the quantum mechanics notation that you'll ever desire.

Suppose that one has two distinguishable particles $A$ and $B$.  One can then write a joint pure-state wave function $|\Psi\rangle$, which is expressed as an expansion in its states weighted by the probability that they occur.  Note that we have written the wave function using Dirac (bra-ket) notation.  It is a column vector, and its Hermitian conjugate is the row vector $\langle \Psi |$.  The prefactor for each term in the expansion of $|\Psi\rangle$ consists of the complex-valued components $C_{ij}$ of an $m \times n$ probability matrix $C$, which satisfies $\mbox{tr}(CC^\dagger) = \mbox{tr}(C^\dagger C) = 1$.  (Recall that $X^\dagger$ refers to the Hermitian conjugate of the matrix $X$.)

Applying the SVD of $C$ (i.e., letting $C = USV^\dagger$, where $U$ and $V$ are unitary matrices\footnote{A unitary matrix $U$ satisfies $UU^\dagger = 1$ and is the complex-valued generalization of an orthogonal matrix.}) and transforming to a single-particle basis allows one to diagonalize $|\Psi\rangle$, which is said to be {\it entangled} if more than one singular value is nonzero.  One can even measure the entanglement using the two-particle density matrix $\rho := |\Psi\rangle \langle\Psi |$ that is given by the outer product of the wave function with itself.  One can then compute the von Neumann entanglement entropy
\begin{equation}
	\sigma = - \sum_{k = 1}^{\mbox{min}(n,m)}|S_k^2|\ln |S_k^2|\,.
\end{equation}
Because $|S_k^2| \in [0,1]$, the entropy is zero for unentangled states and has the value $\ln[\mbox{min}(n,m)]$ for maximally entangled states.

The SVD is also important in other aspects of quantum information.  For example, it can be used to help construct measurements that are optimized to distinguish between a set of (possibly nonorthogonal) quantum states \cite{eldar}.  

%%%%%%%%%%%%%%	

%------------------------

\section{Can You Take Me Higher?}

\begin{center}
\begin{figure}[h!]
\centering
\fbox{\phantom{$\begin{array}{c}x\\x\\x\\x\end{array}$}
\epsfig{file=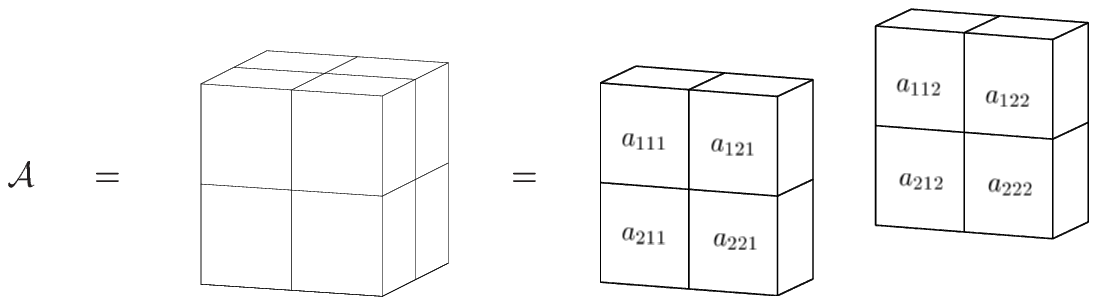}}
\caption{Illustration of a $2\times 2\times 2$ tensor as a cube of data.
This figure originally appeared in Ref.~\cite{KilmerMartin} and is used with permission from Elsevier.
}
\label{fig:tensorpic}
\end{figure}
\end{center}

As we have discussed, the SVD permeates numerous applications and is vital to data analysis.  Moreover, with the availability of cheap memory and advances in instrumentation and technology, it is now possible to collect and store enormous quantities of data for science, medical, and engineering applications.  A byproduct of this wealth is an ever-increasing abundance of data that is fundamentally three-dimensional or higher.  The information is thus stored in multiway arrays---i.e., as tensors---instead of as matrices. An order-$p$ tensor $\mathcal{A}$ is a multiway array with $p$ indices: 
\[\mathcal{A}=(a_{i_1i_2\dots i_p})\in\mathbb{R}^{n_1\times n_2\times\dots\times n_p}\,.\]
Thus, a first-order tensor is a vector, a second-order tensor is a matrix, a third-order tensor is 
a ``cube", and so on.  See Fig.~\ref{fig:tensorpic} for an illustration of a $2\times 2\times 2$
tensor.

Applications involving operations with tensors are now widespread.  They include chemometrics \cite{app1}, psychometrics \cite{Kroon}, signal processing \cite{Comon,Lath_signalproc,Sidiropoulos}, computer vision \cite{tensorface3,tensorface4,tensorface1}, data mining \cite{datamine1,Savas}, networks \cite{Kolda4,multislice}, neuroscience \cite{neuro1,neuro2,neuro3}, and many more.  For example, the
facial recognition algorithm \emph{Eigenfaces} \cite{Turk1,Turk2,eigenface} has been extended to \emph{TensorFaces} \cite{tensorface3}.  To give another example, experiments have shown that fluorescence (i.e., the emission of light from a substance) is modeled well using tensors, as the data follow a trilinear model \cite{app1}.

A common thread in these applications is the need to manipulate the data, usually by compression,  
by taking advantage of its multidimensional structure (see, for example, the recent article \cite{OST}).  Collapsing multiway data to matrices and using standard linear algebra to answer questions about the data often has undesirable consequences.  It is thus important to consider the multiway data directly.

Here we provide a brief overview of two types of higher-order extensions of the matrix SVD.
For more information, see the extensive article on tensor
decompositions \cite{SIREV} and references therein.  Recall from \eqref{svd2} that the SVD is a rank-revealing decomposition.  
The outer product $u_iv_i^T$ in equation \eqref{svd2}
is often written using the notation $u_i\circ v_i$.  Just as
the outer product of two vectors is a rank-1 matrix, the outer product
of three vectors is a rank-1 third-order tensor.  For example, if
$x\in\mathbb{R}^{n_1}$, $y\in\mathbb{R}^{n_2}$, and $z\in\mathbb{R}^{n_3}$, then the outer product $x\circ y\circ z$ has dimension $n_1\times n_2\times n_3$ and is a rank-1 third-order tensor whose $(i,j,k)$th entry is given by $x_iy_jz_k$.  Likewise, an outer product of four vectors
gives a rank-1 fourth-order tensor, etc.  For the rest of this discussion, we will limit our
exposition to third-order tensors, but the concepts generalize easily to order-$p$ tensors.

The \emph{tensor rank} $r$ of an order-$p$ tensor $\mathcal{A}$ is the minimum number of 
rank-1 tensors that are needed to express the tensor.  For a third-order tensor $\mathcal{A}\in\mathbb{R}^{n_1\times n_2\times n_3}$, this implies the representation
\begin{align}
  \label{PARAFAC}
  	\mathcal{A} = \sum_{i=1}^{r} \sigma_i(u_i \circ v_i \circ w_i)\,,
\end{align}
where $\sigma_i$ is a scaling constant.  The scaling constants are the nonzero elements of an
$r\times r\times r$ diagonal tensor $S=(\sigma_{ijk})$. (As discussed in Ref.~\cite{SIREV}, a tensor is called \emph{diagonal} if the only nonzero entries occur in elements $\sigma_{ijk}$ with
$i=j=k$.)  The vectors $u_i$, $v_i$, and $w_i$ are the $i$th columns from
matrices $U\in\mathbb{R}^{n_1\times r}$, $V\in\mathbb{R}^{n_2\times r}$, and $W\in\mathbb{R}^{n_3\times r}$, respectively.  

One can think of equation \eqref{PARAFAC} as an extension of the matrix SVD.  Note, however, 
the following differences.
\begin{enumerate}
\item{The matrices $U$, $V$, and $W$ in (\ref{PARAFAC}) are {\em not} constrained to be orthogonal. Furthermore, an orthogonal decomposition of the form \eqref{PARAFAC} does not exist,
except in very special cases \cite{denis}.
}
\item{The maximum possible rank of a tensor is not given directly from
the dimensions, as is the case with matrices.\footnote{The maximum possible rank
of an $n_1\times n_2$ matrix is $\min(n_1,n_2)$.}  However, loose upper bounds on rank do exist for higher-order tensors.  Specifically, the maximum
possible rank of an $n_1\times n_2\times n_3$ tensor is bounded by 
$\min(n_1n_2,n_1n_3,n_2n_3)$ in general \cite{multiway} and $\floor{3n/2}$ in the case
of $n \times n \times 2$ tensors \cite{jaja,multiway,martinrank,tenberge}.  In
practice, however, the rank is typically much less than these upper bounds. For example,
Ref.~\cite{Comon2} conjectures that the rank of a particular $9\times 9\times 9$ tensor is 19 or 20.
}
\item{Recall that the best rank-$k$ approximation to a matrix is given by the $k$th
partial sum in the SVD expansion (Theorem \ref{thm-EY}). 
However, this result does not extend to higher-order 
tensors.  In fact, the best rank-$k$ approximation to a tensor might not even exist \cite{degen1,degen2}.
}
\item{There is no known closed-form solution to determine the rank $r$
of a tensor \emph{a priori}; in fact, the problem is NP-hard \cite{NPhard}.
Rank determination of a tensor is a widely-studied problem \cite{SIREV}.  
}
\end{enumerate}

In light of these major differences, there exists more than one 
higher-order version of the
matrix SVD.  The different available decompositions are motivated by the application areas.
A decomposition of the form (\ref{PARAFAC}) is called a 
CANDECOMP-PARAFAC (CP) decomposition 
(CANonical DECOMPosition or PARAllel FACtors
model) \cite{CANDECOMP,PARAFACMODEL}, whether or not $r$ is known to be minimal.  However, since an orthogonal decomposition of the form
\eqref{PARAFAC} does not always exist, a \emph{Tucker3} form is often used to guarantee the existence of an orthogonal decomposition as well as 
to better model certain data \cite{NagyKilmer,Savas,tensorface3,tensorface1,
tensorface4}.  

If $\mathcal{A}$ is an $n_1 \times n_2\times n_3$ tensor, then its \emph{Tucker3 decomposition} has the form \cite{Tucker3}
\begin{align}
\label{TUCKER}
  	\mathcal{A} = \sum_{i=1}^{m_1}\sum_{j=1}^{m_2}\sum_{k=1}^{m_3} \sigma_{ijk}
(u_i \circ v_j \circ w_k)\,,
\end{align}
where $u_i$, $v_j$, and $w_k$ are the $i$th, $j$th, and $k$th columns of 
the matrices $U\in\mathbb{R}^{n_1\times m_1}$, $V\in\mathbb{R}^{n_2\times m_2}$, and $W\in\mathbb{R}^{n_3\times m_3}$.  Often, $U$, $V$, and $W$ have orthonormal columns.  The
tensor $S=(\sigma_{ijk})\in\mathbb{R}^{m_1\times m_2\times m_3}$
is called the \emph{core tensor}.  In general, the core tensor $S$ is dense and the decomposition
(\ref{TUCKER}) does not reveal its rank.  Equation \eqref{TUCKER} has also been called the higher-order SVD (HOSVD) \cite{DDV}, though the term ``HOSVD" actually refers 
to a method for computation \cite{SIREV}.  Reference \cite{DDV} demonstrates that the HOSVD is a convincing extension of the matrix SVD.  The HOSVD is guaranteed to exist, and it computes \eqref{TUCKER} directly by calculating the SVDs of the three matrices obtained by ``flattening" the tensor into matrices in each dimension and then using those results to assemble the core tensor.  Yet another extension of the matrix SVD factors a tensor as a product of tensors rather than as an outer product of vectors \cite{KilmerMartin,Martin2}.

%%%%%%%%%%%%%

\section{Everywhere You Go, Always Take the SVD With You.}

The SVD is a fascinating, fundamental object of study that has provided a great deal of of insight into a diverse array of problems, which range from social network analysis and quantum information theory to applications in geology.  The matrix SVD has also served as the foundation from which to conduct data analysis of multiway data by using its higher-dimensional tensor versions.  The abundance of workshops, conference talks, and journal papers in the past decade on multilinear algebra and tensors also demonstrates the explosive growth of applications for tensors and tensor SVDs.  
%Recent attention by the linear algebra community has advanced both the theory and computational aspects of higher-order tensors, but much
%work remains to be done.  
The SVD is an omnipresent factorization in a plethora of application areas.  We recommend it highly.

%%%%%%%%%%%%%%

\paragraph{Acknowledgements.}
We thank Roddy Amenta, Keith Briggs, Keith Hannabuss, Peter Mucha, Steve Simon, Gil Strang, Nick Trefethen, Catalin Turc, and Charlie Van Loan for useful discussions and comments on drafts of this paper.  We also thank Mark Newman for assistance with Figs.~\ref{senate} and \ref{votesvd}.

%%%%%%%%%%%%%%%%%

\bigskip

\noindent\textbf{Carla D. Martin} is an associate professor of mathematics at James Madison University. She has a strong bond with linear algebra and especially with the SVD, inspired by her thesis advisor at Cornell University, Charles Van Loan. She performs research in multilinear algebra and tensors but pretty much loves anything having to do with matrices. She has been a huge promoter of publicizing mathematical applications in industry both in talks and in print, and also serves as the VP of Programs in BIG SIGMAA. She is also an active violinist and enjoys teaching math and music to her three children.
\noindent\textit{Department of Mathematics and Statistics, James
Madison University, Harrisonburg, VA  22807\\carla.dee@gmail.com}

\bigskip

\noindent\textbf{Mason A. Porter} is a University Lecturer in the Mathematical Institute at the University of Oxford.  He is a member of the Oxford Centre for Industrial and Applied Mathematics (OCIAM) and of the CABDyN Complexity Centre.  He is also a Tutorial Fellow of Somerville College in Oxford.  Mason's primary research areas are nonlinear and complex systems, but at some level he is interested in just about everything.  Mason greatly prefers using the word ``myriad" as an adjective rather than as a noun (and it is perfectly correct to do so), but sometimes he must bow to editorial demands.  Mason originally saw the SVD as an undergraduate at Caltech, although it was a particular incident involving an SVD question on a midterm exam in Charlie Van Loan's matrix computations class at Cornell University that helped inspire the original version of this article.    As an exercise (which is related to the solution that he submitted for that exam), he encourages diligent readers to look up the many existing backronyms for SVD.
\noindent\textit{Mathematical Institute, University of Oxford, Oxford, OX1 3LB, United Kingdom\\ porterm@maths.ox.ac.uk}


\begin{thebibliography}{9}

\bibitem{datamine1}
{ E.~Acar, S.~A. \c{C}amtepe, M.~S. Krishnamoorthy, B.~Yener}, 
  Modeling and multiway analysis of chatroom tensors, in 
  \emph{Intelligence and Security Informatics, Lecture Notes in Computer Science}, Vol. 3495, Edited by Kantor, Paul and Muresan, Gheorghe and Roberts, Fred and Zeng, Daniel and Wang, Fei-Yue and Chen, Hsinchun and Merkle, Ralph, Springer, Berlin/Heidelberg, 2005. 181--199.

\bibitem{genome2}
{ O.~Alter, P.~O. Brown, D.~Botstein}, {Singular value decomposition
  for genome-wide expression data processing and modeling}, \emph{Proceedings of the
  National Academy of Sciences} \textbf{97} (2000) 10101--10106.

\bibitem{Amenta2}
{ R.~Amenta, A.~Ewing, A.~Jensen, S.~Roberts, K.~Stevens, M.~Summa,
  S.~Weaver, P.~Wertz}, {A modeling approach to understanding the role
  of microstructure development on crystal-size distributions and on recovering
  crystal-size distributions from thin slices}, \emph{American Mineralogist} \textbf{92} (2007) 1936--1945.

\bibitem{Amenta1}
{ R.~Amenta, B.~Wihlem}, Application of singular value decomposition
  to estimating grain sizes for crystal size distribution analysis,
\emph{GAC-MAC-SEG-SGA Ottawa 2011}, available at \url{http://www.gacmacottawa2011.ca/}.

\bibitem{LAPACK}
{ E.~Anderson, Z.~Bai, C.~Bischof, S.~Blackford, J.~Demmel, J.~Dongarra,
  J.~Du~Croz, A.~Greenbaum, S.~Hammarling, A.~McKenney, D.~Sorensen}, \emph{LAPACK User's Guide}, third edition. Society for Industrial and Applied Mathematics, Philadelphia, 1999.

\bibitem{neuro1}
{ C.~Beckmann, S.~Smith}, {Tensorial extensions of the independent
  component analysis for multisubject {f}{M}{R}{I} analysis}, 
  \emph{NeuroImage} \textbf{25}
  (2005) 294--311.

\bibitem{robotics}
{ C.~Belta, V.~Kumar}, {An {SVD}-based projection method for
  interpolation on ${SE}(3)$}, \emph{IEEE Transactions on Robotics and Automation} \textbf{18}
  (2002) 334--345.

\bibitem{khachiyan}
{ A.~Bengt, R.~E. Stone}, {Khachiyan's linear programming algorithm},
  \emph{Journal of Algorithms} \textbf{1} (1980) 1--13.

\bibitem{tenberge}
{ J.~M.~F. ten Berge}, {Kruskal's polynomial for $2\times 2\times 2$
  arrays and a generalization to $2\times n\times n$ arrays}, \emph{Psychometrika} \textbf{56} (1991) 631--636.

\bibitem{Berry92largescale}
{ M.~W. Berry}, {Large scale sparse singular value computations},
  \emph{International Journal of Supercomputer Applications} \textbf{6} (1992) 13--49.

\bibitem{text}
{ M.~W. Berry, S.~T. Dumais, G.~W. O'Brien}, {Using linear algebra
  for intelligent information retrieval}, \emph{SIAM Review} 
  \textbf{37} (1995) 573--595.

\bibitem{org}
{ J.~R. Boyce, D.~P. Bischak}, {The role of political parties in the
  organization of {C}ongress}, \emph{The Journal of Law, Economics, \& Organization}
  \textbf{18} (2002) 1--38.

\bibitem{CANDECOMP}
{ J.~D. Carroll, J.~Chang}, {Analysis of individual differences in
  multidimensional scaling via an {$N$}-way generalization of
  ``{E}ckart-{Y}oung'' decomposition}, \emph{Psychometrika} 
  \textbf{35} (1970) 283--319.

\bibitem{dhillonsvd}
{ A.~K. Cline, I.~S. Dhillon}, {Computation of the singular value
  decomposition}, in \emph{Handbook of Linear Algebra}, Edited by L.~Hogben, CRC Press, Boca Raton, FL, 2006. 45.1--45.13.

\bibitem{Comon}
{ P.~Comon}, {Tensor decompositions}, in \emph{Mathematics in Signal
  Processing {V}}, Edited by J.~G. McWhirter and I.~K. Proudler, 
  Clarendon Press, Oxford, 2002. 1--24.

\bibitem{Comon2}
{ P.~Comon, J.~M.~F. ten Berge, L.~De Lathauwer, J.~Castaing},
  Generic and typical ranks of multi-way arrays, \emph{Linear Algebra and its
  Applications} \textbf{430} (2009) 2997--3007.

\bibitem{Lath_signalproc}
{ L.~De~Lathauwer, B.~De~Moor}, {From matrix to tensor: Multilinear
  algebra and signal processing}, in \emph{Mathematics in Signal Processing IV}, Edited by
  J.~McWhirter and I.K.~Proudler, Clarendon Press, Oxford, 1998. 1--15.

\bibitem{DDV}
{ L.~De~Lathauwer, B.~De~Moor, J.~Vandewalle}, {A multilinear
  singular value decomposition}, \emph{SIAM Journal of Matrix Analysis and
  Applications} \textbf{21} (2000) 1253--1278.


\bibitem{degen1}
{ V.~De~Silva, L.-H. Lim}, Tensor rank and the ill-posedness of the
  best low-rank approximation problem, \emph{SIAM Journal on Matrix Analysis and
  Applications} \textbf{30} (2008) 1084--1127.

\bibitem{denis}
{ J.~B. Denis, T.~Dhorne}, {Orthogonal tensor decomposition of 3-way
  tables}, in \emph{Multiway Data Analysis}, Edited by R.~Coppi and S.~Bolasco, Elsevier, Amsterdam, 1989. 31--37.

\bibitem{sharad}
{ P.~Diaconis, S.~Goel, S.~Holmes}, {Horseshoes in multidimensional
  scaling and local kernel methods}, \emph{Annals of Applied Statistics} \textbf{2} (2008) 777--807.

\bibitem{dieci}
{ L.~Dieci, M.~G. Gasparo, A.~Papini}, {Path following by {SVD}}, in
  \emph{Computational Science - ICCS 2006, Lecture Notes in Computer Science}, Vol. 3994, 2006. 677--684.

\bibitem{EY}
{ G.~Eckart, G.~Young}, {The approximation of one matrix by another
  of lower rank}, \emph{Psychometrika} \textbf{1} (1936) 211--218.

\bibitem{eldar}
{ Y.~C. Eldar, G.~D. {Forney, Jr.}}, {On quantum detection and the
  square-root measurement}, \emph{IEEE Transactions on Information Theory} \textbf{47} (2001) 858--872.

\bibitem{fennpca}
{ D.~J. Fenn, M.~A. Porter, S.~Williams, M.~McDonald, N.~F. Johnson,
  N.~S. Jones}, {Temporal evolution of financial market correlations},
  \emph{Physical Review E} \textbf{84} (2011) 026109.

\bibitem{GolubReinsch}
{ G.~H. Golub, C.~Reinsch}, { Singular value decomposition and least
  squares solutions}, \emph{Numerische Mathematik} \textbf{14} (1970) 403--420.

\bibitem{vanloan}
{ G.~H. Golub, C.~F. Van~Loan}, \emph{Matrix Computations}, third edition, The Johns
  Hopkins University Press, Baltimore, MD, 1996.

\bibitem{PARAFACMODEL}
{ R.~A. Harshman}, {Foundations of the {PARAFAC} procedure: Model and
  conditions for an `explanatory' multi-mode factor analysis}, 
  \emph{UCLA Working
  Papers in phonetics} \textbf{16} (1970) 1--84.

\bibitem{NPhard}
{ J.~Hastad}, {Tensor rank is {NP}-complete}, \emph{Journal of Algorithms} \textbf{11}
  (1990) 664--654.

\bibitem{higgins06}
{ M.~D. Higgins}, \emph{Quantitative Textural Measurements in Igneous and
  Metamorphic Petrology}, Cambridge University Press, Cambridge, U.K., 2006.

\bibitem{genome1}
{ N.~S. Holter, M.~Mitra, A.~Maritan, M.~Cieplak, J.~R. Banavar, N.~V.
  Fedoroff}, {Fundamental patterns underlying gene expression profiles:
  Simplicity from complexity}, \emph{Proceedings of the National Academy of Sciences}
 \textbf{ 97} (2000) 8409--8414.

\bibitem{jaja}
{ J.~Ja'Ja'}, {Optimal evaluation of pairs of bilinear forms}, \emph{SIAM
  Journal on Computing} \textbf{8} (1979) 443--461.

\bibitem{KilmerMartin}
{ M.~E. Kilmer, C.~D. Martin}, {Factorization strategies for
  third-order tensors}, \emph{Linear Algebra and its Applications} 
  \textbf{435} (2011) 641--658.

\bibitem{chemphys}
{ T.~Kinoshita}, {Singular value decomposition approach for the
  approximate coupled-cluster method}, \emph{Journal of Chemical Physics} \textbf{119} (2003) 7756--7762.

\bibitem{Kolda4}
{ T.~G. Kolda, B.~W. Bader}, Higher-order web link analysis using
  multilinear algebra, in {\em Data Mining {ICDM 2005},} Proceedings of the 5th IEEE International Conference, IEEE Computer Society, 2005. 242--249.

\bibitem{SIREV}
\leavevmode\vrule height 2pt depth -1.6pt width 23pt, Tensor
  decompositions and applications, \emph{SIAM Review} \textbf{51} (2009) 455--500.

\bibitem{Kroon}
{ P.~M. Kroonenberg}, \emph{Three-Mode Principal Component Analysis: Theory
  and Applications}, DSWO Press, Leiden, 1983.

\bibitem{multiway}
{ J.~B. Kruskal}, {Rank, decomposition, and uniqueness for 3-way and
  {$N$}-way arrays}, in \emph{Multiway Data Analysis}, Edited by
  R.~Coppi and S.~Bolasco, Elsevier, Amsterdam, 1989. 7--18.

\bibitem{Larsen98lanczosbidiagonalization}
{ R.~M. Larsen}, {Lanczos bidiagonalization with partial
  reorthogonalization}, in {\em Efficient Algorithms for Helioseismic Inversion}, Ph.D. Thesis, Department of Computer Science, University of Aarhus, 1998. Part II, Chapter A.

\bibitem{Rokhlin}
{ E.~Liberty, F.~Woolfe, P.-G. Martinsson, V.~Rokhlin, M.~Tygert},
  Randomized algorithms for the low-rank approximation of matrices,
  \emph{Proceedings of the National Academy of Sciences} \textbf{104} (2007)
  20167--20172.

\bibitem{susanbook}
{ K.~Marida, J.~T. Kent, S.~Holmes}, \emph{Multivariate Analysis},
  Academic Press, New York, 2005.

\bibitem{marsh98}
{ A.~D. Marsh}, On the interpretation of crystal size distributions in
  magmatic systems, {\em Journal of Petrology} \textbf{39} (1998)  553--599.

\bibitem{martinrank}
{ C.~D. Martin}, The rank of a $2\times 2\times 2$ tensor, 
\emph{Linear and
  Multilinear Algebra}, \textbf{59} (2011) 943--950.

\bibitem{Martin2}
{ C.~D. Martin, R.~Shafer, B.~LaRue}, A recursive idea for
  multiplying order-$p$ tensors, 
  %{\em SIAM Journal on Scientific Computing}, 
  submitted (2011).

\bibitem{neuro2}
{ E.~Mart\'{i}nez-Montes, P.~A. Vald\'{e}s-Sosa, F.~Miwakeichi, R.~I.
  Goldman, M.~S. Cohen}, Concurrent {EEG/fMRI} analysis by multiway
  partial least squares, {\em NeuroImage} \textbf{22} (2004) 1023--1034.

\bibitem{neuro3}
{ F.~Miwakeichi, E.~Mart\'{i}nez-Montes, P.~A. Vald\'{e}s-Sosa,
  N.~Nishiyama, H.~Mizuhara, Y.~Yamaguchi}, Decomposing {EEG} data
  into space-time-frequency components using parallel factor analysis,
  {\em NeuroImage} \textbf{22} (2004) 1035--1045.

\bibitem{moshtagh}
{ N.~Moshtagh}, Minimum volume enclosing ellipsoid-- From MATLAB Central
  File Exchange (1996) MathWorks, \url{http://www.mathworks.com/matlabcentral/fileexchange/9542-minimum-volume-enclosing-ellipsoid}.

\bibitem{multislice}
{ P.~J. Mucha, T.~Richardson, K.~Macon, M.~A. Porter, J.-P. Onnela},
  Community structure in time-dependent, multiscale, and multiplex
  networks, {\em Science} \textbf{328} (2010) 876--878.

\bibitem{eigenface}
{ N.~Muller, L.~Magaia, B.~M. Herbst}, Singular value
  decomposition, eigenfaces, and 3{D} reconstructions, {\em SIAM Review} \textbf{46} (2004) 518--545.

\bibitem{NagyKilmer}
{ J.~G. Nagy, M.~E. Kilmer}, Kronecker product approximation for
  preconditioning in three-dimensional imaging applications, {\em IEEE Trans. Image
  Proc.} \textbf{15} (2006) 604--613.

\bibitem{OST}
{ I.~V. Oseledets, D.~V. Savostianov, E.~E. Tyrtyshnikov}, Tucker
  dimensionality reduction of three-dimensional arrays in linear time, {\em SIAM
  Journal of Matrix Analysis and Applications} \textbf{30} (2008)  939--956.

\bibitem{degen2}
{ P.~Paatero}, Construction and analysis of degenerate {PARAFAC}
  models, {\em Journal of Chemometrics}, \textbf{14} (2000) 285--299.

\bibitem{pasyou}
{ R.~Pa\v{s}kauskas, L.~You}, Quantum correlations in two-boson wave
  functions, {\em Physical Review A} \textbf{64} (2001) 042310.

\bibitem{voteview}
{ K.~T. Poole}, Voteview, Department of Political Science, University of Georgia, 2011, \url{http://voteview.com}.

\bibitem{pr97}
{ K.~T. Poole, H.~Rosenthal}, \emph{Congress: A Political-Economic History
  of Roll Call Voting}, Oxford University Press, Oxford, 1997.

\bibitem{conglong}
{ M.~A. Porter, P.~J. Mucha, M.~E.~J. Newman, A.~J. Friend}, 
  Community structure in the {U}nited {S}tates {H}ouse of {R}epresentatives,
  {\em Physica A} \textbf{386} (2007) 414--438.

\bibitem{congshort}
{ M.~A. Porter, P.~J. Mucha, M.~E.~J. Newman, C.~M. Warmbrand}, 
A  network analysis of committees in the {U}nited {S}tates {H}ouse of
  {R}epresentatives, {\em Proceedings of the National Academy of Sciences} \textbf{102}
  (2005) 7057--7062.

\bibitem{preskill}
{ J.~Preskill}, Lecture notes for physics 229: Quantum information and
  computation (2004), California Institute of Technology, available at 
  \url{http://www.theory.caltech.edu/~preskill/ph229/}.

\bibitem{Amenta4}
{ A.~Randolph, M.~Larson}, \emph{Theory of Particulate Processes},
  Academic Press, New York, 1971.

\bibitem{saari}
{ D.~G. Saari}, ed., \emph{Decisions and Elections: Explaining the
  Unexpected}, Cambridge University Press, Cambridge, 2001.

\bibitem{Savas}
{ B.~Savas, L.~Eld\'{e}n}, Handwritten digit classification using
  higher-order singular value decomposition, {\em Pattern Recognition} 
  \textbf{40} (2007) 993--1003.

\bibitem{sch}
{ J.~Schliemann, D.~Loss, A.~H. MacDonald}, Double-occupancy
  errors, adiabaticity, and entanglement of spin qubits in quantum dots,
  {\em Physical Review B} \textbf{63} (2001) 085311.

\bibitem{Sidiropoulos}
{ N.~Sidiropoulos, R.~Bro, G.~Giannakis}, Parallel factor analysis
  in sensor array processing, {\em IEEE Transactions on Signal Processing} \textbf{48} (2000) 2377--2388.

\bibitem{app1}
{ A.~Smilde, R.~Bro, P.~Geladi}, \emph{Multi-way Analysis: Applications
  in the Chemical Sciences}, Wiley, 2004.

\bibitem{Stewart93}
{ G.~W. Stewart}, On the early history of the singular value
  decomposition, {\em SIAM Review} \textbf{35} (1993) 551--566.

\bibitem{strang}
{ G.~Strang}, \emph{Linear Algebra and its Applications}, fourth edition, Thomson Brooks/Cole, 2005.

\bibitem{Tucker3}
{ L.~R. Tucker}, Some mathematical notes on three-mode factor
  analysis, {\em Psychometrika} \textbf{31} (1966) 279--311.

\bibitem{Turk1}
{ M.~Turk, A.~Pentland}, Eigenfaces for recognition, {\em Journal of
  Cognitive Neuroscience} \textbf{3} (1991a).

\bibitem{Turk2}
\leavevmode\vrule height 2pt depth -1.6pt width 23pt, Face recognition
  using {E}igenfaces, {\em Proc. of Computer Vision and Pattern Recognition} \textbf{3}
  (1991b) 586--591.

\bibitem{twain}
{ M.~Twain}, Pudd'nhead {W}ilson's new calendar, in {\em Following the
  Equator}, Samuel L. Clemens, Hartford, CT, 1897.

\bibitem{tensorface3}
{ M.~A.~O. Vasilescu, D.~Terzopoulos}, Multilinear analysis of image
  ensembles: Tensorfaces, {\em Computer Vision -- ECCV 2002}, Proceedings of the 7th European Conference, Lecture Notes in Computer Science, Vol. 2350, 2002. 447--460.

\bibitem{tensorface1}
\leavevmode\vrule height 2pt depth -1.6pt width 23pt, Multilinear image
  analysis for face recognition, {\em Pattern Recognition--ICPR 2002, Proceedings of the International Conference}, Vol. 2, 2002. 511--514.

\bibitem{tensorface4}
\leavevmode\vrule height 2pt depth -1.6pt width 23pt, Multilinear subspace
  analysis of image ensembles, {\em Computer Vision and Pattern Recognition--CVPR 2003}, Proceedings of the 2003 IEEE Computer Society
  Conference, 2003. 93--99.

\bibitem{compress}
{ L.~Xu, Q.~Liang}, Computation of the singular value
  decomposition, in \emph{Wireless Algorithms, Systems, and Applications 2010,
  Lecture Notes in Computer Science}, Vol. 6221, Edited by 
  G.~Pandurangan, V.~S.~A.
  Kumar, G.~Ming, Y.~Liu, and Y.~Li, Springer-Verlag, Berlin, 2010.
  338--342.

\end{thebibliography}
\end{document}